\documentclass[12pt]{amsart} 
\usepackage{algorithm} 
\usepackage{algorithmic}
\usepackage[text={6.5in,9in},centering]{geometry}

\usepackage{amsmath} 
\usepackage{amssymb} 
\usepackage{amsfonts}
\usepackage{amscd} 
\usepackage{color}
\usepackage{graphicx}

\newcommand{\GG}{\mathcal{G}} 
\newcommand{\VV}{\mathcal{V}}
\newcommand{\calS}{\mathcal{S}}

\newcommand{\EE}{\mathcal{E}}

\def\argmax{\arg\max}
\def\argmin{\arg\max}

\usepackage{subfig}
\usepackage{bm} 
\usepackage{bbm}
\newcommand{\bone}{\ensuremath{\mathbbm{1}}}
\newcommand{\Reals}[1]{{\rm I\! R}^{#1}}

\newtheorem{theorem}{Theorem}[section]
\newtheorem{lemma}[theorem]{Lemma}

\newtheorem{example}[theorem]{Example}

\numberwithin{equation}{section} \numberwithin{table}{section}
\begin{document}

 \title{Parallel Unsmoothed Aggregation Algebraic Multigrid Algorithms
   on GPUs}
\thanks{The first, third and fourth author were 
 supported in part by the grants NSF DMS-1217142 and OCI-0749202 and DoE DE-SC0006903}

  \author{James Brannick}
\address{Department of Mathematics, The
      Pennsylvania State University, University Park, PA 16802, USA.}
 \email{brannick@psu.edu}
\author{Yao Chen}
\address{Microsoft Corporation, 1065 La Avenida, Mountain View,  CA
  94043, USA.}
\email{yaoc@microsoft.com}
\author{Xiaozhe Hu}
\address{Department of Mathematics, The
      Pennsylvania State University, University Park, PA 16802, USA.}
\email{hu\_x@math.psu.edu} 
\author{Ludmil  Zikatanov}
\address{Department of Mathematics, The
      Pennsylvania State University, University Park, PA 16802, USA.}
\email{ludmil@psu.edu} 

\begin{abstract}
We design and implement a parallel algebraic
  multigrid method for isotropic graph Laplacian problems on multicore
  Graphical Processing Units (GPUs).  The proposed AMG
  method is based on the aggregation framework.  The setup phase of
  the algorithm uses a parallel maximal independent set algorithm in
  forming aggregates and the resulting coarse level hierarchy is then
  used in a K-cycle iteration solve phase with a $\ell^1$-Jacobi
  smoother.  Numerical tests of a parallel implementation of the
  method for graphics processors are presented to demonstrate its
  effectiveness.  
\end{abstract}

\maketitle

\section{Introduction} We consider development of a multilevel
iterative solver for large-scale sparse linear systems corresponding
to graph Laplacian problems for graphs with balanced vertex degrees. A
typical example is furnished by the matrices corresponding to the
(finite difference)/(finite volume)/(finite element) discretizations of
scalar elliptic equation with mildly varying coefficients on unstructured grids.

Multigrid (MG) methods have been shown to be very efficient iterative
solvers for graph Laplacian problems and numerous parallel MG
solvers have been developed for such systems.  Our aim here is to
design an algebraic multigrid (AMG) method for solving the graph Laplacian system
and discuss the implementation of such methods on
multi-processor parallel architectures, with an emphasis on
implementation on Graphical Processing Units (GPUs).

The programming environment which we use in this paper is the Compute
Unified Device Architecture (CUDA) toolkit introduced in 2006 by
NVIDIA which provides a framework for programming on GPUs. Using this
framework in the last 5 years several variants of Geometric Multigrid
(GMG) methods have been implemented on GPUs
\cite{Goodnight2003,Bolz2003,Goddeke2008,Grossauer2008,Feng2010,2012FengC_ShuS_XuJ_ZhangC-aa}
and a high level of parallel performance for the GMG algorithms on
CUDA-enabled GPUs has been demonstrated in these works.

On the other hand, designing AMG methods for massively parallel
heterogenous computing platforms, e.g., for clusters of GPUs, is very
challenging mainly due to the sequential nature of the coarsening
processes (setup phase) used in AMG methods.  In most AMG algorithms,
coarse-grid points or basis are selected sequentially using graph
theoretical tools (such as maximal independent sets and graph
partitioning algorithms).  Although extensive research has been
devoted to improving the performance of parallel coarsening
algorithms, leading to notable improvements on CPU architectures
\cite{Cleary1998,Tuminaro2000,Krechel2001,Henson2002,DeSterck2006,Chow2006,Joubert2006,Tuminaro2000}, on a single GPU \cite{Bell2011,Kraus, Haase2010}, and on multiple GPUs \cite{Emans.M;Liebmann.M;Basara.B.2012a}, the setup phase is still
considered a bottleneck in parallel AMG methods.  We mention the
work in~\cite{Bell2011}, where a smoothed aggregation setup is
developed in CUDA for GPUs.

In this paper, we describe a parallel AMG method based on the
un-smoothed aggregation AMG (UA-AMG) method.  The setup algorithm we
develop and implement has several notable design features.  A key
feature of our parallel aggregation algorithm is that it first chooses
coarse vertices using a parallel maximal independent set algorithm
\cite{DeSterck2006} and then forms aggregates by grouping coarse
level vertices with their neighboring fine level vertices, which, in turn, avoids
ambiguity in choosing fine level vertices to form aggregates.  Such a design
eliminates both the memory write conflicts and conforms to the CUDA
programming model.  The triple matrix product needed to compute the
coarse-level matrix (a main bottleneck in parallel AMG setup
algorithms) simplifies significantly in the UA-AMG setting, reducing
to summations of entries in the matrix on the finer level.  The parallel
reduction sums available in CUDA are quite an efficient tool for this
task during the AMG setup phase.  Additionally, the UA-AMG setup
typically leads to low grid and operator complexities.

In the solve phase of the proposed algorithm, a
K-cycle~\cite{panayot-book,amli-1,amli-2} is used to accelerate the 
convergence rate of the multilevel UA-AMG method.  Such multilevel 
method optimizes the coarse grid correction and results in
an approximate two-level method.  Two parallel relaxation schemes
considered in our AMG implementation are a damped Jacobi smoother and
a parameter free $\ell^1$-Jacobi smoother
 introduced in~\cite{2009KolevT_VassilevskiP-aa}  and its weighted
 version in~\cite{2012BrezinaM_VassilevskiP-aa}. 
To further accelerate the convergence rate of the resulting K-cycle
method we apply it as a preconditioner to a nonlinear conjugate
gradient method.


The remainder of the paper is organized as follows.  In Section
\ref{sec:UA-AMG}, we review the UA-AMG method.  Then, in Section
\ref{sec:setup}, a parallel graph aggregation method is introduced,
which is our main contribution.  The parallelization of the solve
phase is discussed in Section \ref{sec:solve}. In Section
\ref{sec:numerics}, we present some numerical results to demonstrate
the efficiency of the parallel UA-AMG method.  

\section{Unsmoothed Aggregation AMG} \label{sec:UA-AMG} 

The linear system of interest has as coefficient matrix the graph
Laplacian corresponding to an undirected connected graph
$\GG=(\VV,\EE)$. Here, $\VV$ denotes the set of vertices and $\EE$
denotes the set of edges of $\GG$.  We set $n=|\VV|$ (cardinality of
$\VV$). By $(\cdot,\cdot)$ we denote the inner product in
$\ell^2(\Reals{n})$ and the superscript $t$ denotes the adjoint with
respect to this inner product.  The {\em graph Laplacian}
$A:\Reals{n}\mapsto \Reals{n}$ is then defined via the following
bilinear form
\begin{eqnarray*}
(A u, v)=\sum_{k=(i,j)\in\EE} \omega_{ij}(u_{i}-u_{j})(v_{i}-v_{j}) +
\sum_{j\in\calS} \omega^D_j u_jv_j, \quad\calS\subset \VV. 
\end{eqnarray*}
We assume that the weights $\omega_{ij}$, and $\omega^D_j$ are
strictly positive for all $i$ and $j$. The first summation is over the
set of edges $\EE$ (over $k\in \EE$ connecting the vertices $i$ and
$j$), and
$u_{i}$ and $u_{j}$ are the $i$-th and $j$-th coordinate of the vector
$u\in \Reals{n}$, respectively. We also assume that the subset of vertices $\calS$ is such
that the resulting matrix $A$ is symmetric positive definite (SPD). If
the graph is connected $\calS$ could contain only one vertex and $A$
will be SPD. For matrices corresponding to the discretization scalar
elliptic equation on unstructured grids $\calS$ is the set of vertices
near (one edge away from) the boundary of the computational domain.  The linear system of
interest is then
\begin{equation}\label{eqn:linear system}
A u = f.
\end{equation}

With this system of equation we associate a multilevel hierarchy which
consists of spaces $V_0\subset V_1\subset\ldots\subset V_L=\Reals{n}$,
each of the spaces is defined as the range of
interpolation/prolongation operator $P_{l-1}^l:\Reals{n_{l-1}}\mapsto
V_{l}$ with $\operatorname{Range}(P_{l-1}^{l})=V_{l-1}$. 

Given the $l$-th level matrix $A_l \in \mathbb{R}^{n_l  \times n_l}$, the aggregation-based prolongation matrix
$P_{l-1}^{l}$ is defined in terms of a non-overlapping partition of
the $n_l$ unknowns at level $l$ into the $n_{l-1}$ nonempty disjoint
sets $G_j^{l}$, $j=1,\dots, n_{l-1}$, called aggregates.  An algorithm
for choosing such aggregates is presented in the next section. The
prolongation $P_{l-1}^{l}$ is the $n_{l} \times n_{l-1}$ matrix with
columns defined by partitioning the constant vector, $\bone=(1,\ldots,1)^t$, with
respect to the aggregates:
\begin{equation} \label{def:P} (P_{l-1}^{l})_{ij} =
\begin{cases} 1 & \text{if} \ i \in G_j^{l} \\ 0 & \text{otherwise}
\end{cases} \quad i = 1, \dots, n_l, \quad j = 1, \dots, n_{l-1}.
\end{equation} The resulting coarse-level matrix
$A_{l-1} \in \mathbb{R}^{n_{l-1} \times n_{l-1}}$ is then defined by
the so called ``triple matrix product'', namely, 
\begin{equation} \label{def:Ac} A_{l-1} = (P_{l-1}^{l})^t A_l
(P_{l-1}^{l}).
\end{equation} 
Note that since we consider UA-AMG, the interpolation operators are
boolean matrices such that the entries in the coarse-grid matrix $A_{l-1}$ can
be obtained from a simple summation process:
\begin{equation}\label{eqn:acsum} (A_{l-1})_{ij} = \sum_{s \in G_{i}}
\sum_{t \in G_{j}} a_{st}, \quad i,j = 1, 2, \cdots, n_{l-1}.
\end{equation} 
Thus, the triple matrix product, typically \emph{the} costly procedure in an AMG setup,
simplifies significantly for UA-AMG to reduction sums.

We now introduce a general UA-AMG method (see Algorithm \ref{alg:UA-AMG})
and in the subsequent sections we describe the implementation 
of each of the components of Algorithm~\ref{alg:UA-AMG} for GPUs.

\begin{algorithm}[htbp]
\caption{UA-AMG} \label{alg:UA-AMG} \flushleft{{\bf Setup Phase:}}
\begin{algorithmic}[1] \STATE Given $n_0$ (size of the coarsest level) and $L$ (maximum levels) 
\STATE $l \gets L$, 
\WHILE { $N_{l} \geq n_0 \ \& \ l
> 0$ } \STATE Construct the aggregation
$\mathcal{N}_{l}^{i}$, $i=1,2, \dots, N_{l+1}$ based on $A_{l}$,
\STATE Compute $A_{l-1}$ by~\eqref{eqn:acsum}, \STATE $l \gets l-1$,
\ENDWHILE
\end{algorithmic} {\bf Solve Phase:}
\begin{algorithmic}[1] \IF { (On the coarsest level) } \STATE solve
$A_{l}u_{l} = f_{l}$ exactly, 
\ELSE 
\STATE Pre-smoothing: $u_{l} \gets
\texttt{smooth}(u_{l}, A_{l}, f_{l})$, 
\STATE Restriction: compute
$r_{l-1} = (P_{l-1}^l)^T (f_l - A_l u_l)$, 
\STATE Coarse grid
correction: solve $A_{l-1}e_{l-1} = r_{l-1}$ approximately by recursively
calling the AMG on coarser level $l-1$ and get $e_{l-1}$, 
\STATE
Prolongation: compute $u_{l} \gets u_{l} + P_{l-1}^l e_{l-1}$, 
\STATE
Post-smoothing: $u_{l} \gets \texttt{smooth}(u_{l}, A_{l}, f_{l})$.
\ENDIF
\end{algorithmic}\end{algorithm}

\section{The Setup Phase} \label{sec:setup}

Consider the system of linear equations~\eqref{eqn:linear system}
corresponding to an unweighted graph $\GG=\{\VV,\EE\}$ 
partitioned into two subgraphs $\GG_{k}=\{\VV_{k},\EE_{k}\}, k=1,2$.
Further assume that the two subgraphs are stored on separate computes.
To implement a Jacobi or Gauss-Seidel smoother for the graph Laplacian
equation with respect to $\GG$, the communication between the two
computers is proportional to the number of edge cuts of such a
partitioning, given by
$$
|\EE \backslash (\EE_{1}\cup\EE_{2})| .
$$
Therefore, a partition corresponding to the minimal edge cut in the
graph results in the fastest implementation of such smoothers.  
This in turn gives a heuristic argument, as also
suggested in \cite{metis1}, \cite{metis2}, that when partitioning the graph
in subgraphs (aggregates) the subgraphs should have a
similar number of vertices and have a small ``perimeter.''  Such a
partitioning can be constructed by choosing any vertex in the graph,
naming it as a coarse vertex, and then aggregating it with its
neighboring vertices.  This heuristic motivates our
aggregation method. The
algorithm consists of a sequence of two subroutines: first, a parallel
maximal independent set algorithm is applied to identify coarse
vertices; then a parallel graph aggregation algorithm follows, so that
subgraphs (aggregates) centered at the coarse vertices are formed.

In the algorithm, to reduce repeated global memory read access and write conflicts, we
impose explicit manual scheduling on data caching and flow control in
the implementations of both algorithms; the aim is to achieve the
following goals:
\begin{enumerate}
\item (Read access coalescence): To store the data that a node uses
frequently locally or on a fast connecting neighboring node.
\item (Write conflicting avoidance): To reduce, or eliminate the
situation that several nodes need to communicate with a center node
simultaneously.
\end{enumerate}

\subsection{A Maximal Independent Set Algorithm}

The idea behind such algorithm is to simplify the memory coalescence,
and design a random aggregation algorithm where there are as many as
possible threads loading from a same memory location, while as few as
possible threads writing to a same memory location.  Therefore, it is
natural to have one vertex per thread when choosing the coarse
vertices.  For vertices that are connected the corresponding
processing threads should be wrapped together in a group.  By doing so, repeated
memory loads from the global memory can be avoided.

However, we also need to ensure that no two coarse vertices compete
for a fine level point, because either atomic operations as well as inter-thread
communication is costly on a GPU.  Therefore, the coarse vertices are
chosen in a way that any two of them are of distance 3 or more, which
is the same as finding a maximal independent set of vertices for the graph
corresponding to $A^{2}$, where $A$ is the graph Laplacian of a given
graph $\GG$, so that each fine level vertex can be determined independently
which coarse vertex it associates with.

Given an undirected unweighted graph $\GG=\{\VV,\EE\}$, we first find
a set $C$ of coarse vertices such that
\begin{eqnarray}\label{d>=3} d(i,j) \geq 3 ,\qquad \forall i,\,j\in C,
i\neq j .
\end{eqnarray} 
Here, $d(\cdot,\cdot)$ is the graph distance function defined
recursively as
\begin{eqnarray*} d(i,j) = \left\{
\begin{matrix} 0 , & i=j ; \\ \displaystyle\min_{k:(i,k)\in E}
d(k,j)+1 , & i\neq j .
\end{matrix} \right.
\end{eqnarray*} 
Assume we obtain such set $C$, or even a subset of $C$, we can then
form aggregates, by picking up a vertex $i$ in $C$ and defining an
aggregate as a set containing $i$ and its neighbors.  The condition
\eqref{d>=3} guarantees that two distinct vertices in $C$, do not
share any neighbors.  The operation of marking the
numberings of subgraphs on the fine grid vertices is write conflict
free, and the restriction imposed by \eqref{d>=3} ensures that
aggregates can be formed independently, and simultaneously.  

The rationale of the independent set algorithm is as follows: 
First, a random vector $v$ is generated, each component of which
corresponds to a vertex in the graph.  Then we define the set
$C$ as the following:
\begin{eqnarray*} C = \big\{ i \mid v_{i}>v_{j}, \forall j :
0<d(i,j)<3 \big\} .
\end{eqnarray*} 
If $C$ is not empty, then such construction results in a collection of
vertices in $C$ is of distance 3 or more.  Indeed, assume that
$d(i,j)<3$ for $i,\,j\in C$, let $v_{i}>v_{j}$.  From the definition
of the set $C$, we immediately conclude that $i\notin C$. Of course,
more caution is needed when $C$ defined above is empty (a situation
that may occur depending on the vector $v$). However, this can be
remedied, by assuming that the vector $v$ (with random entries) has a
global maximum, which is also a local maximum. The $C$ contains at
least this vertex. The same algorithm can be applied then recursively
to the remaining graph (after this vertex is removed). In practice,
$C$ does not contain one but more vertices.

\subsection{Parallel Graph Aggregation Algorithm}

We here give a description of the parallel aggregation algorithm,
running the exact copies of the code on each thread.

\begin{algorithm}[h!]
\caption{Parallel Aggregation Algorithm (PAA)} \label{alg:paa}
\begin{enumerate}
\item[(1)] Generate a quasi-random number and store it in $v_{i}$, as
\[ v_{i} \longleftarrow \text{quasi\_random}(i) ;
\] mark vertex $i$ as ``unprocessed''; wait until all threads complete
these operations.
\item[(2)]
\begin{enumerate}
\item[(2a)] Goto (2d) if $i$ is marked ``processed'', otherwise
continue to (2b).
\item[(2b)] Determine if the vertex $i$ is a coarse vertex, by check
if the following is true.
\[ v_{i} > v_{j} , \qquad \forall j : (A^{2})_{ij}\neq 0 \text{ and
$j$ is unprocessed } .
\] If so, continue to (2c); if not, goto (2d).
\item[(2c)] Form an aggregate centered at $i$.  Let $S_{i}$ be a set
of vertices defined as
\[ S_{i} =\big\{ j \mid v_{i} \geq v_{j} , \forall j :
(A^{2})_{ij}\neq 0 \text{ and $j$ is unprocessed } \big\} .
\] Define a column vector $w$ such that
\[ w_{k} = \left\{
\begin{matrix} 1, & k \in S_{i} ; \\ 0, & k \not \in S_{i} .
\end{matrix} \right.
\] Mark vertices $j\in S_{i}$ ``processed'' and request an atomic
operation to update the prolongator $P$ as
\[ P \longleftarrow [P, w] \;.
\]
\item[(2d)] Synchronize all threads (meaning: wait until all
threads reach this step).
\item[(2e)] Stop if $i$ is marked ``processed'', otherwise goto step
  (2a).
\end{enumerate}
\end{enumerate}
\end{algorithm}

Within each pass of the Parallel Aggregation Algorithm (PAA, Algorithm \ref{alg:paa}), 
the following two steps are applied to each vertex $i$.
\begin{enumerate}
\item[(A)] Construct a set $C$ which contains coarse vertices.
\item[(B)] Construct an aggregate for each vertex in $C$.
\end{enumerate} 
Note that these two subroutines can be executed in a parallel
fashion. Indeed, step (A) does not need to be applied to the whole
graph before starting step (B).  Even if $C$ is partially completed,
any operation in step (B) will not interfere step (A), running on the
neighboring vertices and completing the construction of $C$.  A
problem for this approach is that it usually cannot give a set of
aggregates that cover the vertex set $V$ after 1 pass of step (A) and
step (B).  We thus run several passes and the algorithm terminates
when a complete cover is obtained.  The number of
passes is reduced if we make the set $C$ as large as possible in each
pass, therefore the quasi-random vector $v$ needs to have a lot of
local maximums.  Another heuristic argument is that $C$ needs to be
constructed in a way that every coarse vertex has a large number of
neighboring vertices.  Numerical experiments suggest that the following is
a good way of generating the vector $v$ with the desired properties.
\begin{eqnarray}\label{quasi-random} v_{i} \longleftarrow
\text{quasi\_random}(i) := d_{i}+\big( (i \text{ mod } 12) +
\text{rand}()\big)/12 .
\end{eqnarray} where $d_{i}$ is the degree of the vertex $i$, and
$\text{rand}()$ generates a random number uniformly distributed on the
interval $[0,1]$.

\subsection{Aggregation Quality Improvements}
\label{ssec:rank1}

To improve the quality of the aggregates, we can either impose some
constrains during the aggregation procedure (which we call in-line
optimization), or introduce a post-process an existing aggregation in
order to improve it. One in-line strategy that we use to improve the
quality of the aggregation is to limit the number of vertices in an
aggregate during the aggregation procedure.  However, such limitations
may result in a small coarsening ratio. In such case, and numerical
results suggest that applying aggregation process twice, which is equivalent to
in skipping a level in a multilevel hierarchy, can compensate that.

Our focus is on a post-processing strategy,  which we name
``rank one optimization''. It uses an \emph{a priori} estimate to
adjust the interface (boundary) of a pair of aggregates, so that the
aggregation based two level method, with a fixed smoother, converges
fast locally on those two aggregates.

We consider the connected graph formed by a union of aggregates (say a
pair of them, which will be the case of interest later), and let
$\widehat{n}$ be the dimension of the underlying vector space.  Let
$\widehat{A}:\Reals{\widehat{n}}\mapsto\Reals{\widehat{n}}$ be a
semidefinite weighted graph Laplacian (representing a local
sub-problem) and $\widehat R$ be a given local smoother.  As is usual
for semidefinite graph Laplacians, we consider the subspace
$\ell^2$-orthogonal to the null space of $\widehat{A}$ and we denote
it by $V$. The $\ell^2$ orthogonal projection on $V$ is denoted here
by $\Pi_V$. Let $\widehat S=I-\widehat R \widehat A$ be the error
propagation operator for the smoother $\widehat{R}$.  We consider the
two level method whose error propagation matrix is
\[
E(V_c)=E(V_c; \widehat{S})=
 (I-Q_{\widehat A}(V_c))(I-\widehat R\widehat A).
\]
Here $V_c\subset V$ is a subspace and $Q_{\widehat A}(V_c)$ is the
$\widehat A$-orthogonal projection of the elements of $V$ onto the
coarse space $V_c$.  In what follows we use the notation
$E(V_{c};\widehat{S})$ when we want to emphasize the dependence on
$\widehat{S}$. We note that $Q_{\widehat{A}}(V_c)$ is well defined on
$V$ because $\widehat{A}$ is SPD on $V$ and hence it
$(\widehat{A}\cdot,\cdot)$ is an inner product on $V$. We also have
that $Q_{\widehat{A}}(V_c)$ self-adjoint on $V$ and under the
assumption $V_c\subset V$, we obtain
$Q_{\widehat{A}}(V_c)=\Pi_VQ_{\widehat{A}}(V_c)$ and
$\Pi_VQ_{\widehat{A}}(V_c)=Q_{\widehat{A}}(V_c)\Pi_V$. 
Also, $\widehat{S}_V=\Pi_V\widehat{S}$ is self-adjoint on $V$ in the
$(\widehat{A}\cdot,\cdot)$ inner product iff $\widehat{R}$ is
self-adjoint in the $\ell^2$-inner product on $\Reals{\widehat{n}}$.

  We now introduce the operator $T(V_c)$ (recall that $V_c\subset V$)
\[
T(V_c)= T(V_c;\widehat{S})= 
\widehat S -
E(V_c)=Q_{\widehat{A}}(V_c)(I-\widehat{R}\widehat{A}) = Q_{\widehat A}(V_c)\widehat{S}.
\]
and from the definition of  $Q_{\widehat A}$ for all $v\in V$ we have
\begin{equation}\label{e=s-t}
  |E(V_c)v|_{\widehat A}^{2}=|\Pi_VE(V_c)v|_{\widehat A}^{2}=
|\Pi_V\widehat Sv|_{\widehat A}^{2}-|\Pi_VT(V_c)v|_{\widehat A}^{2}
=|\widehat{S}_Vv|_{\widehat A}^{2}-|T(V_c)v|_{\widehat A}^{2}.
\end{equation}
We note the
following identities which follow directly from the definitions above
and the assumption $V_c\subset V$:
\begin{equation}\label{eq:reduction}
|E(V_c;\widehat{S})|_{\widehat{A}}=|\Pi_VE(V_c;\widehat{S}_V)|_{\widehat{A}},
\quad
|T(V_c;\widehat{S})|_{\widehat{A}}=|\Pi_VT(V_c;\widehat{S}_V)|_{\widehat{A}}.
\end{equation}

The relation \eqref{e=s-t} suggests that, in order to minimize the
seminorm $|E(V_c)v|_{\widehat A}$ with respect to the coarse space
$V_c$, we need to make $|T(V_c)|_{\widehat{A}}$ maximal. The following
lemma quantifies this observation and is instrumental in showing how
to optimize locally the convergence rate when the subspaces $V_c$ are
one dimensional.  In the statement of the lemma we use $\arg\min$ to
denote a subset of minimizers of a given, not necessarily linear,
functional $F(x)$ on a space $X$. More precisely, we set
\[
y\in \arg\min_{x\in X}{F(x)},\quad\mbox{if and only if},\quad 
F(y)=\min_{x\in X}{F(x)}.
\]
We have similar definition (with obvious changes) for the set
$\displaystyle\argmin_{x\in X} F(x)$.

\begin{lemma}\label{the-only-lemma}
  Let $\widehat{S}_V=\Pi_V\widehat{S}$, be the projection of the local
  smoother on $V$,  and $\mathcal V_c$ be the set of all one dimensional
  subspaces of $V$. Then we have  the following:
\begin{eqnarray}
&&|\widehat{S}_V|_{\widehat{A}}=\max_{V_c \in \mathcal V_c}|T(V_c)|_{\widehat{A}}, \label{1d-thm-1.1} \\
&&\mbox{If $W_c\in \displaystyle\argmax_{V_c \in \mathcal{V}_c}|T(V_c)|_{\widehat{A}}$, then
$W_c\in \displaystyle\arg\min_{V_c\in \mathcal V_c}|E(V_c)|_{\widehat{A}}$}, \label{1d-thm-1.2} 
\end{eqnarray}
 where $E(V_c)= (I-Q_{\widehat{A}}(V_c))\widehat S$ and  $T(V_c)=Q_{\widehat{A}}(V_c)\widehat S$. 
\end{lemma}

\noindent\textbf{Proof.} 
From the identities~\eqref{eq:reduction} it follows that we can
restrict our considerations on $V\subset\Reals{\widehat{n}}$ and that
we only need to prove the Lemma with  $E(V_c)=\Pi_VE(V_c;\widehat{S}_V)$ 
and $T(V_c)=\Pi_VT(V_c;\widehat{S}_V)$.
In order to make the presentation more transparent, we denote
$|\cdot|=|\cdot|_{\widehat{A}}$, $\Pi=Q_{\widehat{A}}$. 
Let us mention also that by orthogonality
in this proof we mean orthogonality in the $(\widehat{A}\cdot,\cdot)$ inner
product on $V$.   The proof then proceeds as follows. 

Let $\varphi\in V$ be such that 
$|\widehat{S}_V\varphi|=|\widehat{S}_V||\varphi|$. We set
$W_c=\operatorname{span}\{\widehat{S}_V\varphi\}$. 
Note that for such choice of $W_c$ we have 
$\Pi(W_c)\widehat{S}_V \varphi=\widehat{S}_V \varphi$ and hence
\[
|\widehat{S}_V|=\frac{|\widehat{S}_V\varphi|}{|\varphi|} =
\frac{|T(W_c)\varphi|}{|\varphi|}\le |T(W_c)| .
\]
On the other hand, for all $V_c\in\mathcal{V}_c$ we have
$|\Pi(V_c)|=1$ and we then conclude that
\begin{equation}\label{eq:all-inequalities}
|T(V_c)|=|\Pi(V_c)\widehat{S}_V|\leq|\Pi(V_c)||\widehat{S}_V|= |\widehat{S}_V|\le |T(W_c)|.
\end{equation}
By taking a maximum on $\mathcal{V}_c$ in \eqref{eq:all-inequalities}, 
we conclude the following thus prove \eqref{1d-thm-1.1}. 
\[
|T(W_{c})| 
\leq \max_{V_c \in \mathcal{V}_c}|T(V_c)|
\leq |\widehat{S}_V|\le |T(W_c)|.
\]

To prove \eqref{1d-thm-1.2}, we observe that for
any $W_c \in \displaystyle\argmax_{V_c \in\mathcal{V}_c}|T(V_c)|$ 
the inequalities in~\eqref{eq:all-inequalities} become equalities and hence 
\[
|\Pi(W_c)\widehat{S}_V|=|\widehat{S}_V|=|\widehat{S}_V\Pi(W_c)| .
\]
This implies that $\displaystyle|\widehat{S}_V|=\max_{w\in  W_c}\frac{|\widehat{S}_Vw|}{|w|}$. It is also clear that $|\widehat{S}_Vw|=|\widehat{S}_V||w|$
for all $w\in W_c$, 
because $W_c$ is one dimensional. In addition, since $\widehat{S}_V$ is
self-adjoint, it follows that $W_c$ is the span of the eigenvector of $\widehat{S}_V$
with eigenvalue of magnitude $|\widehat{S}_V|$. 
Next, for any $V_c\in \mathcal{V}_c$ we have
\begin{equation*}
  |E(V_c)|=|(I-\Pi(V_c))\widehat{S}_V|=
|\widehat{S}_V(I-\Pi(V_c))|=\max_{v\in V_c^{\perp}}\frac{|\widehat{S}_Vv|}{|v|}.
\end{equation*}
By the mini-max principle 
(see \cite[pp. 31-35]{1989CourantR_HilbertD-aa} or~\cite{1974HalmosP-aa}) 
we have that  $|E(V_c)|\ge \sigma_2$, 
where $\sigma_2$ is the second largest singular value of $\widehat{S}_V$ and with
equality holding iff $V_c=W_c$. 
This completes the proof. \hfill $\blacksquare$

We now move on to consider a pair of aggregates.  Let $\widehat{A}$ be
the graph Laplacian of a connected positively weighted graph
$\widehat{\GG}$ which is union of two aggregates $\VV_{1}$ and
$\VV_{2}$.  Furthermore, let $\bone_{\VV_1}$ be the characteristic
vector for $\VV_1$, namely a vector with components equal to
$1$ at the vertices of $\VV_1$ and equal to zero at the vertices of
$\VV_2$. Analogously we have a characteristic vector 
$\bone_{\VV_2}$ for
$\VV_2$. Finally, let $V_c(\VV_{1},\VV_{2})$ be the space of
vectors that are linear combinations of $\bone_{\VV_1}$ and 
$\bone_{\VV_2}$.
More specifically, the subspace $V_c$ is defined as
\[
V_c(\VV_{1},\VV_{2})=\text{span}\left\{
\left(\frac{1}{|\VV_1|}\bone_{\VV_{1}}-\frac{1}{|\VV_2|}\bone_{\VV_{2}}\right)\right\}, 
\]
Let $\mathcal V_c$ be the set of subspaces defined above for all
possible pairs of $\VV_{1}$ and $\VV_{2}$, such that
$\widehat{\GG}=\VV_{1}\cup\VV_{2}$. Note that 
by the definition above, every pair $(\VV_1,\VV_2)$ 
gives us a space $V_c \in \mathcal{V}_c$ which is orthogonal to the null space of
$\widehat{A}$, i.e. orthogonal to $\bone=\bone_{\VV_1}+\bone_{\VV_2}$. 

We now apply the result of Lemma~\ref{the-only-lemma} and show how to
improve locally the quality of the partition (the convergence rate
$|E(V_{c})|_{\widehat{A}}$) by reducing the problem of minimizing the
$\widehat{A}$-norm of $E(V_{c})$ to the problem of finding the maximum
of the $\widehat{A}$-norm of the rank one transformation $T(V_{c})$.
Under the assumption that the spaces $V_c$ are orthogonal to the null
space of $\widehat{A}$ (which they satisfy by construction) from
Lemma~\ref{the-only-lemma} we conclude that the spaces $W_c$ which
minimize $|E(V_c)|_{\widehat{A}}$ also maximize
$|T(V_c)|_{\widehat{A}}$.

For the pair of aggregates, $|T(V_c)|_{\widehat A}$ is the largest
eigenvalue of
$\widehat{S}^TAQ_{\widehat{A}}(V_c)\widehat{S}\widehat{A}^{\dagger}$,
where $A^{\dagger}$ is the pseudo inverse of $A$. Clearly, the matrix
$\widehat{S}^TAQ_{\widehat{A}}(V_c)\widehat{S}\widehat{A}^{\dagger}$,
is also a rank one matrix and hence
\[
|T(V_c)|_{\widehat A}=\text{tr}(\widehat{S}^TAQ_{\widehat{A}}(V_c)\widehat{S}\widehat{A}^{\dagger}).
\]

During optimization steps, we calculate the trace using the fact that
for any rank one matrix $W$ we have
\begin{equation}\label{rank1-trace}
\text{tr}(W)=\frac{W_k^{T}W_k}{W_{kk}} = \frac{W_k^{T}W_k}{e_k^{T}W_{k}e_{k}},
\end{equation}
where $W_{kk}$ is a nonzero diagonal entry (any nonzero diagonal
entry) and $W_k$ is the $k$-th column of $W$. The
formula~\eqref{rank1-trace} is straightforward to prove if we set
$W=uv^T$ for two column vectors $u$ and $v$, and also suggests a
numerical algorithm. We devise a loop computing $W_{k}=We_{k}$, and
$W_{kk}=e_{k}^{T}W_k e_{k}$, for $k=1,\dots,m$, where $m$ is the
dimension of $W$; The loop is terminated whenever $W_{kk}\neq 0$, and
we compute the trace via \eqref{rank1-trace} for this $k$. In
particular for the examples we have tested,
$W=\widehat{S}^{T}\widehat{A}Q_{\widehat{A}}(V_c)\widehat{S}\widehat{A}^{\dagger}$
is usually a full matrix and we observed that the loop almost always
terminated when $k=1$.

The algorithm which traverses all pairs of neighboring aggregates and
optimizes their shape is as follows.

\begin{algorithm}[h!]
\caption{Subgraph Reshaping Algorithm} \label{alg:sra}
\begin{enumerate}
\item[]
\textit{Input}: 
Two set of vertices, $\VV_{1}$ and $\VV_{2}$, 
corresponding to a pair of neighboring subgraphs. 
\textit{Output}: 
Two sets of vertices, $\tilde \VV_{1}$ and $\tilde \VV_{2}$ satisfying that 
\[
\tilde \VV_{1} \cup \tilde \VV_{2} = \VV_{1} \cup \VV_{2} 
,\qquad \text{ and } \qquad 
\big| |\tilde \VV_{1}|-|\tilde \VV_{2}| \big| \leq 1, 
\]
and the subgraphs corresponding to $\tilde \VV_{1}$ and $\tilde \VV_{2}$ are both connected. 
\item[(1)]
Let $n=|\tilde \VV_{1}|+|\tilde \VV_{2}|$, then compute 
$m=\lfloor n/2 \rfloor$. 
\item[(2)]
Run in parallel to generate all partitionings such that the vertices set 
\[
\tilde \VV_{1} \cup \tilde \VV_{2} = \VV_{1} \cup \VV_{2} 
, \qquad 
|\tilde \VV_{1}|=m, 
\]
and the subgraphs derived by $\tilde \VV_{1}$ and $\tilde \VV_{2}$ are connected. 
\item[(3)]
Run in parallel to compute the norm $|T(V_c)|_{\widehat A}$ for all partitionings get from step (2), 
and return the partitioning that results in maximal $|T(V_c)|_{\widehat A}$. 
\end{enumerate}
\end{algorithm}

The subgraph reshaping algorithm fits well the programming model of a multicore
GPU.  We demonstrate this algorithm on two example problems, and later
show its potential as a post-process for the parallel aggregation
algorithm (Algorithm~\ref{alg:paa}) outlined in the previous
section. In the examples that follow next we use the rank one
optimization and then measure the quality of the coarse space also by
computing the energy norm of the $|Q|_{\widehat{A}}$, where 
$Q$ is the $\ell^{2}$-orthogonal projection to the space $W_c$. 

\begin{example}\label{example-1}: Consider a graph Laplacian $\widehat A$ corresponding to a
  graph which is a $4\times 4$ square grid.  The weights on the edges
  are all equal to $1$.  We start with an obviously non-optimal
  partitioning as shown on the left of Figure \ref{fig:isotropic}, of
  which the resulting two level method, consisting of $\ell^1$-Jacobi
  pre- and post-smoothers and an exact coarse level solver, has a
  convergence rate $|E|_{\widehat A}=0.84$, and $|Q|_{\widehat
    A}^{2}=1.89$.  After applying Algorithm \ref{alg:sra}, the refined
  aggregates have the shapes shown on the right of Figure
  \ref{fig:isotropic}, of which the two level method has the same
  convergence rate $|E|_{\widehat A}=0.84$ but the square of the
  energy seminorm is reduced to $|Q|_{\widehat A}^{2}=1.50$.

\begin{figure}[h]
\captionsetup{width=30em}
\begin{center}
\parbox{90pt}{\includegraphics[page=1,width=90pt]{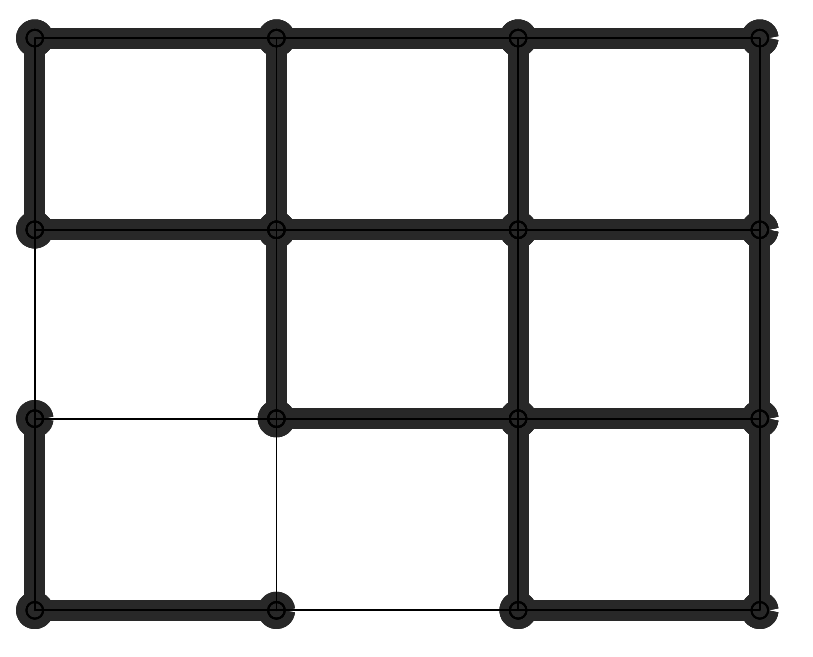}}
$\qquad\longrightarrow\qquad$
\parbox{90pt}{\includegraphics[page=2,width=90pt]{figs_gray}}
\end{center}
\caption{Subgraph reshaping algorithm applied on a graph 
representing an isotropic coefficient elliptic PDE.}\label{fig:isotropic}
\end{figure}
\end{example}

\begin{example}\label{example-2}
Consider a graph Laplacian $\widehat A$ corresponding to
a graph which is a $4\times 4$ square grid, 
on which all horizontal edges are weighted $1$ while all vertical edges are weighted $10$.  
Such graph Laplacian represents anisotropic coefficient elliptic equations with 
Neumann boundary conditions. 
Start with a non-optimal partitioning as shown on the left of Figure \ref{fig:anisotropic}, 
of which the resulting two level method has a convergence rate $|E|_{\widehat A}=0.96$ 
and $|Q|_{\widehat A}^{2}=4.88$.  
After applying Algorithm \ref{alg:sra}, 
the refined aggregates have the shapes shown on the right of Figure \ref{fig:anisotropic}, 
of which the two level convergence rate is reduced to $|E|_{\widehat A}=0.90$ 
and the energy of the coarse level projection is also reduced as $|Q|_{\widehat A}^{2}=1.50$. 

\begin{figure}[h]
\captionsetup{width=30em}
\begin{center}
\parbox{90pt}{\includegraphics[page=1,width=90pt]{figs_gray}}
$\qquad\longrightarrow\qquad$
\parbox{90pt}{\includegraphics[page=3,width=90pt]{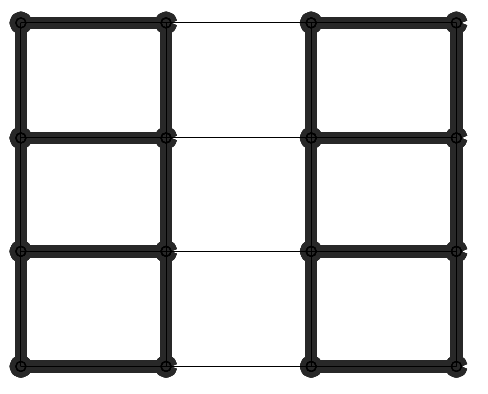}}
\end{center}
\caption{Subgraph reshaping algorithm applied on a graph 
representing an anisotropic coefficient elliptic PDE.}\label{fig:anisotropic}
\end{figure}
\end{example}

\section{Solve Phase} \label{sec:solve} In this section, we discuss
the parallelization of the solver phase on GPU.  More precisely, we
will focus on the parallel smoother, prolongation/restriction, MG
cycle, and sparse matrix-vector multiplication.

\subsection{Parallel Smoother} An efficient parallel smoother is
crucial for the parallel AMG method.  For the sequential AMG method,
Gauss-Seidel relaxation is widely used and has been shown to have a
good smoothing property.  However the standard Gauss-Seidel is a
sequential procedure that does not allow efficient parallel
implementation.  To improve the arithmetic intensity of the smoother
and make it work better with SIMT based GPUs, we adopt the well-known
Jacobi relaxation, and introducing a damping factor to improve the
performance of the Jacobi smoother.  For a matrix $A \in \mathbb{R}^{n
\times n}$ and its diagonals are denoted by $D = \text{diag}(a_{11},
a_{22}, \cdots, a_{nn})$, the Jacobi smoother can be written in the
following matrix form
\begin{equation*} x^{m+1} = x^m + \omega D^{-1} r^m, \quad
\text{where} \ r^m = b-Ax^m,
\end{equation*} or component-wise
\begin{equation*} x^{m+1}_i = x^m_i + \omega a_{ii}^{-1} r^m_i.
\end{equation*} This procedure can be implemented efficiently on GPUs
by assigning one thread to each component, and update the
corresponding components locally and simultaneously.

We also consider the so-called $\ell^1$ Jacobi smoother, which is
parameter free.  
Define 
\[
M = \text{diag}(M_{11}, M_{22}, \cdots,
M_{nn}),
\] 
where $M_{ii} = a_{ii} + d_{ii}$ with $d_{ii} = \sum_{j\neq
i} | a_{ij} |$, and the $\ell^1$ Jacobi has the following matrix form
\begin{equation*} x^{m+1} = x^m + M^{-1} r^m, \quad \text{where} \ r^m
= b-Ax^m,
\end{equation*} or component-wise
\begin{equation*} x^{m+1}_i = x^m_i + M_{ii}^{-1} r^m_i.
\end{equation*} 
In~\cite{2009KolevT_VassilevskiP-aa, 2012BrezinaM_VassilevskiP-aa} 
it has been show that if $A$ is symmetric positive definite, the
smoother is always convergent and has multigrid smoothing properties
comparable to full Gauss-Seidel smoother if $a_{ii} \geq \theta
d_{ii}$ and $\theta$ is bounded away from zero.  Moreover, because its
formula is very similar to the Jacobi smoother, it can also be
implemented efficiently on GPUs by assigning one thread to each
component, and update the corresponding the component locally and
simultaneously.

\subsection{Prolongation and Restriction} For UA-AMG method, the
prolongation and restriction matrices are piecewise constant and
characterize the aggregates. Therefore, we can preform the
prolongation and restriction efficiently in UA-AMG method.  Here, the
output array $\texttt{aggregation}$ (column index of $P$), which
contains the information of aggregates, plays an important rule.
\begin{itemize}
\item {\bf Prolongation:} Let $v^{l-1}\in \mathbb{R}^{n_{l-1}}$, so
that the action $ v^l = P_{l-1}^{l}v^{l-1}$ can be written
component-wise as follows:
\[ (v^l)_i = (P_{l-1}^{l}v^{l-1} )_{i} = (v^{l-1})_{j}, \quad j \in
G^{l-1}_{i}
\] Assign each thread to one element of $v^l$, and the array
\texttt{aggregation} can be used to obtain information about $j \in G^{l-1}_i$,
i.e., $i = \texttt{aggregation}[j]$, so that prolongation can be
efficiently implemented in parallel.

\item {\bf Restriction:} Let $v^{l}\in \mathbb{R}^{n_{l}}$, so that
the action $(P_{l-1}^{l})^Tv^{l}$ can be written component-wise as
follows:
\[ (v^{l-1})_i = ((P_{l-1}^{l})^Tv^{l} )_{i} = \sum_{j \in
G^{l-1}_{i}}(v^{l})_{j}.
\] Therefore, each thread is assigned to an element of $v^{l-1}$, and
the array \texttt{aggregation} can be used to obtain information about $j \in
G^{j-1}_i$, i.e., to find all $j$ such that $\texttt{aggregation}[j] =
i$.  By doing so, the action of restriction can also be implemented in
parallel.
\end{itemize}

\subsection{K-cycle} 
Unfortunately, in general, UA-AMG with V-cycle is
not an optimal algorithm in terms of convergence rate. But on the
other hand, in many cases, UA-AMG using two-grid solver phase gives
optimal convergence rate for graph Laplacian problems. This motivated us to use
other cycles instead of V-cycle to mimic the two-grid algorithm. The
idea is to invest more works on the coarse grid, and make the method become closer to an exact two-level method, then hopefully, the resulting cycle will have optimal convergence rate.

The particular cycle we will discuss here is the so-called K-cycle
(Nonlinear AMLI-cycle) and we refer
to~\cite{panayot-book,amli-1,amli-2} 
for details on its implementation in general.

\subsection{Sparse Matrix-Vector Multiplication on GPUs} As
the K-cycle will be used as a preconditioner for 
Nonlinear Preconditioned Conjugate Gradient (NPCG) method, 
the sparse matrix-vector multiplication (SpMV) has  major contribution
to the computational work involved. An efficient SpMV algorithm
on GPU requires a suitable sparse matrix storage format.  How
different storage formats perform in SpMV is extensively studied in
\cite{Bell2008}.  This study shows that the need for coalesce
accessing of the memory makes ELLPACK (ELL) format 
one of the most efficient sparse matrix storage formats on GPUs 
when each row of the sparse matrix has roughly the
same nonzeros.  In our study, because our main focus is on the
parallel aggregation algorithm and the performance of the UA-AMG
method, we still use the compressed row storage (CSR) format, which
has been widely used for the iterative linear solvers on CPU.
Although this is not an ideal choice for GPU implementation, the
numerical results in the next section already show the efficiency of
our parallel AMG method.  

\section{Numerical Tests} \label{sec:numerics} In this section, we
present numerical tests using the proposed parallel AMG methods.
Whenever possible we compare the results with the CUSP libraries
\cite{Garland2010}. CUSP is an open source C++ library of generic
parallel algorithms for sparse linear algebra and graph computations
on CUDA-enabled GPUs.  All CUSP's algorithms and implementations have
been optimized for GPU by NVIDIA's research group.  To the best of our
knowledge, the parallel AMG method implemented in the CUSP package is
the state-of-the-art AMG method on GPU. We use as test
problems several discretizations of the Laplace equation.

\subsection{Numerical Tests for Parallel Aggregation Algorithm}
Define $Q$, the $\ell^2$ projection on the piece-wise constant
space $\operatorname{Range}(P)$, as the following: 
\begin{eqnarray*} 
Q=P(P^{T}P)^{-1}P^{T} .
\end{eqnarray*} 
We present several tests showing how the energy norm of this
projection changes with respect to different parameters used in the 
parallel aggregation algorithm, since the convergence rate is an increasing
function of $\|Q\|_A$. 

The tests involving $\|Q\|_A$ further suggest two additional features necessary to
get a multigrid hierarchy with predictable results.  First, the sizes
of aggregates need to be limited, and second, the columns of the
prolongator $P$ need to be ordered in a deterministic way, regardless
of the order that aggregates are formed.  The first requirement can be
fulfilled simply by limiting the sizes of the aggregates in each pass
of the parallel aggregation algorithm.  We make the second requirement
more specific.  Let $c_{k}$ to be the index of the coarse vertex of the
$k$-th aggregate.  We require that $c_{k}$ should be an increasing
sequence and then use the $k$-th column of $P$ to record the aggregate
with the coarse vertex numbered $c_{k}$.  This can be done by using a
generalized version of the prefix sum algorithm
\cite{Blelloch90prefixsums}.

We first show in Table~\ref{square-dirichlet} the coarsening ratios
(in the parenthesis in the table) and the energy norms $\|Q\|_{A}^{2}$
of a two grid hierarchy, for a Laplace equation with Dirichlet
boundary conditions on a structured grid containing $n^{2}$ vertices.
The limit on the size of an aggregate is denoted by $t$, which
suggests that, any aggregate can include $t$ vertices or less, which
directly implies that the resulting coarsening ratio is less or equal
to $t$.

\begin{table}[h]
\captionsetup{width=30em}
\center
\begin{tabular}{|r|cccc|}
\hline
& $t=2$ & $t=3$ & $t=4$ & $t=5$ \\
\hline
$n=128$ & (1.99) 1.71 & (2.03) 2.04 & (2.41) 2.46 & (2.97) 3.15 \\
$n=256$ & (1.99) 1.72 & (2.39) 2.57 & (2.96) 2.59 & (2.99) 3.20 \\
$n=512$ & (2.00) 1.72 & (2.01) 2.08 & (2.40) 2.48 & (2.99) 3.22 \\
\hline
\end{tabular}
\caption{(coarsening ratios) and energy norm of 
a two grid hierarchy of a Laplace equation on a uniform grid with Dirichlet boundary conditions.}
\label{square-dirichlet}
\end{table}

For the same aggregations on the graphs that represent Laplace
equations with Neumann boundary conditions, the corresponding
coarsening ratio (in parenthesis) and $|Q|_{A}^{2}$ seminorms with respect
to grid size $n$ and limiting threshold $t$ are shown in Table~\ref{square-neumann}.

\begin{table}[h]
\captionsetup{width=30em}
\center
\begin{tabular}{|r|cccc|}
\hline
& $t=2$ & $t=3$ & $t=4$ & $t=5$ \\
\hline
$n=128$ & (1.99) 1.87 & (2.03) 2.11 & (2.41) 2.48 & (2.97) 3.24 \\
$n=256$ & (1.99) 1.74 & (2.39) 2.59 & (2.96) 2.62 & (2.99) 3.24 \\
$n=512$ & (2.00) 1.87 & (2.01) 2.11 & (2.40) 2.49 & (2.99) 3.24 \\
\hline
\end{tabular}
\caption{(coarsening ratios) and energy norm of 
a two grid hierarchy of a Laplace equation on a uniform grid with Neumann boundary conditions.}
\label{square-neumann}
\end{table}

In Table~\ref{half-of-12} we present the computed bounds on the
coarsening ratio and energy of a two level hierarchy\thanks{By energy
  of a two level hierarchy here, we mean the semi-norm of the $\ell^2$
  projection on the coarse space.} when the fine level is an $n\times
n$. Such results are valid for any structured grid with $n^2$ vertices
(not just $n=126,\ldots, 132$) This is seen as follows: (1) From
equation \eqref{quasi-random}, it follows that if we consider two
grids of sizes $n_{1}\times n_{1}$ and $n_{2}\times n_{2}$
respectively, and such that
\[
(n_{1}-n_{2}) \equiv 0 \mod 12
,\quad \text{ or } \quad
(n_{1}+n_{2}) \equiv 0 \mod 12
\]
then our aggregation algorithm results in the same pattern of $C$ points
on these two grids; (2) As a consequence grids of size $n\times n$ for
$n=126\equiv 6 \mod 12$ to $n=132\equiv 0 \mod 12$ give \emph{all}
possible coarsening patterns that can be obtained by our aggregation
algorithm on \emph{any} 2D tensor product grid. As a
conclusion, the values of the coarsening ratios and the energy
semi-norm given in Table~\ref{half-of-12} are valid for any 2D
structured grid. 

\begin{table}[h]
\captionsetup{width=30em}
\center
\begin{tabular}{|r|cccc|}
\hline
& $t=2$ & $t=3$ & $t=4$ & $t=5$ \\
\hline
$n=126$ & (2.00) 2.11 & (2.00) 2.07 & (2.36) 2.73 & (2.36) 2.73 \\
$n=127$ & (1.99) 1.86 & (2.01) 1.98 & (2.01) 2.49 & (2.01) 2.34 \\
$n=128$ & (1.99) 1.71 & (2.03) 2.04 & (2.41) 2.47 & (2.97) 3.15 \\
$n=129$ & (1.99) 1.84 & (2.02) 2.04 & (2.03) 2.31 & (2.02) 2.42 \\
$n=130$ & (1.99) 1.77 & (2.40) 2.21 & (2.92) 2.86 & (2.94) 2.94 \\
$n=131$ & (1.99) 2.61 & (2.01) 2.41 & (2.01) 2.45 & (2.00) 2.49 \\
$n=132$ & (1.98) 2.09 & (2.21) 2.81 & (2.33) 2.89 & (2.26) 2.94 \\
\hline
\end{tabular}
\caption{(coarsening ratios) and energy norm of 
a two grid hierarchy of a Laplace equation on a uniform grid with Neumann boundary conditions.}
\label{half-of-12}
\end{table}

We also apply this aggregation method on graphs corresponding to
Laplace equations on 2 dimensional unstructured grids with Dirichlet
or Neumann boundary conditions.  The unstructured grids are
constructed by perturbing nodes in an $n\times n$ square lattice
($n=128,256,512$), followed by triangulating the set of perturbed
points using a Delaunay triangulation.  The condition numbers of the
Laplacians, derived using finite element discretization of the Laplace
equations on the mentioned unstructured grids with Dirichlet boundary
conditions, are about $1.2\times 10^{4}$, $5.0\times 10^{4}$, and
$2.1\times 10^{5}$ respectively.  The coarsening ratios and
$|Q|_{A}^{2}$ are listed in table \ref{unstructured-dirichlet} and
\ref{unstructured-neumann}.  We remark here
that 
we also apply the parallel aggregation algorithm without imposing
limit on the size of an aggregate, and the corresponding numerical
results are listed in columns named ``$t=\infty$''.

\begin{table}[h]
\captionsetup{width=34em}
\center
\begin{tabular}{|r|ccccc|}
\hline
& $t=2$ & $t=3$ & $t=4$ & $t=5$ & $t=\infty$\\
\hline
$n=128$ & 
(1.80) 2.39 & (2.53) 2.44 & (3.17) 3.10 & (3.73) 3.46 & (4.91) 3.30\\
$n=256$ & 
(1.79) 2.39 & (2.52) 2.60 & (3.15) 3.18 & (3.69) 3.46 & (4.91) 3.41\\
$n=512$ & 
(1.80) 2.38 & (2.55) 2.67 & (3.19) 3.26 & (3.72) 3.56 & (4.93) 3.40\\
\hline
\end{tabular}
\caption{(coarsening ratios) and energy norm of 
a two grid hierarchy of a Laplace equation with Dirichlet boundary conditions 
discretized on an unstructured grid.}
\label{unstructured-dirichlet}
\end{table}

\begin{table}[h]
\captionsetup{width=34em}
\center
\begin{tabular}{|r|ccccc|}
\hline
& $t=2$ & $t=3$ & $t=4$ & $t=5$ & $t=\infty$\\
\hline
$n=128$ & 
(1.80) 2.39 & (2.53) 2.54 & (3.17) 3.18 & (3.73) 3.48 & (4.91) 3.33\\
$n=256$ & 
(1.79) 2.49 & (2.52) 2.65 & (3.15) 3.20 & (3.69) 3.48 & (4.91) 3.41\\
$n=512$ & 
(1.80) 2.47 & (2.55) 2.80 & (3.19) 3.27 & (3.72) 3.57 & (4.93) 3.53\\
\hline
\end{tabular}
\caption{(coarsening ratios) and energy norm of 
a two grid hierarchy of a Laplace equation with Neumann boundary conditions 
discretized on an unstructured grid.}
\label{unstructured-neumann}
\end{table}

We note that in Table \ref{unstructured-dirichlet} and Table \ref{unstructured-neumann}, 
the coarsening ratios are not large enough to result in small operator
complexity. We then estimate the energy norm $|Q|_{A}^{2}$ 
when $Q$ is the $\ell^{2}$ orthogonal projection from any level of the multigrid hierarchy 
to any succeeding sub-levels. 
We start with a Laplace equation on a $128^{2}$ structured square grid, 
set the limit of sizes of aggregates as $t=5$ on each iteration of aggregation, 
and stop when the coarsest level is of less than 100 degrees of freedom. 
If denoting the finest level by level 0, and the coarsest level we get is level 5. 
The coarsening ratios and energy norm $|Q|_{A}^{2}$ between any two levels are shown 
in table \ref{square-dirichlet-allqa} and \ref{square-neumann-allqa}. 

\begin{table}[h]
\captionsetup{width=34em}
\center
\begin{tabular}{|r|ccccc|}
\hline
& 0 & 1 & 2 & 3 & 4\\
\hline
1 & (2.97) 3.15 & - & - & - & - \\
2 & (11.3) 7.34 & (3.81) 3.09 & - & - & - \\
3 & (38.6) 15.3 & (13.0) 6.74 & (3.41) 2.68 & - & - \\
4 & (113.8) 31.9 & (38.3) 13.9 & (10.1) 5.25 & (2.95) 3.91 & -\\
5 & (321.3) 54.5 & (108.2) 22.5 & (28.4) 9.09 & (8.33) 4.53 & (2.82) 4.77 \\
\hline
\end{tabular}
\caption{(coarsening ratios) and energy norms squares of 
a multigrid hierarchy of a Laplace equation with Dirichlet boundary conditions 
discretized on a $128^{2}$ grid.}
\label{square-dirichlet-allqa}
\end{table}

\begin{table}[h!]
\captionsetup{width=34em}
\center
\begin{tabular}{|r|ccccc|}
\hline
& 0 & 1 & 2 & 3 & 4\\
\hline
1 & (2.97) 3.23 & - & - & - & - \\
2 & (11.3) 7.68 & (3.81) 3.28 & - & - & - \\
3 & (38.6) 16.9 & (13.0) 7.99 & (3.41) 2.94 & - & - \\
4 & (113.8) 42.5 & (38.3) 20.6 & (10.1) 7.07 & (2.95) 4.23 & -\\
5 & (321.3) 98.6 & (108.2) 48.3 & (28.4) 17.0 & (8.33) 7.11 & (2.82) 5.75 \\
\hline
\end{tabular}
\caption{(coarsening ratios) and energy norm of 
a multigrid hierarchy of a Laplace equation with Neumann boundary conditions 
discretized on a $128^{2}$ grid.}
\label{square-neumann-allqa}
\end{table}

We observe that, 
on the diagonal of the table \ref{square-dirichlet-allqa} and \ref{square-neumann-allqa}, 
the energy norms are comparable to the coarsening ratios, 
until the last level where the grid becomes highly unstructured. 
This suggests that a linear or nonlinear AMLI solving cycle can give 
both a good convergence rate and a favorable complexity. 
It is also observed that, 
on the lower triangular part of the table \ref{square-dirichlet-allqa}-\ref{square-neumann-allqa}, 
the energy norms are always smaller than the corresponding coarsening ratios, 
which motives us to design in the future flexible cycles that detect and skip unnecessary levels. 

\begin{figure}[h!]
\captionsetup{width=30em}
\begin{center}
\parbox{220pt}{\includegraphics[page=4,width=220pt]{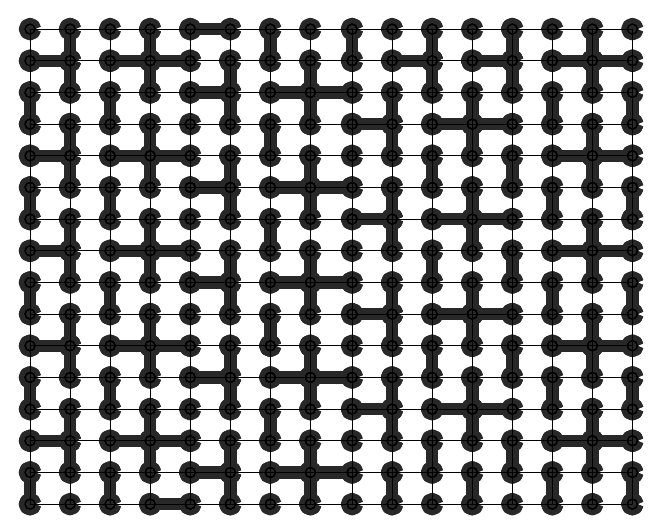}}
$\longrightarrow$
\parbox{220pt}{\includegraphics[page=5,width=220pt]{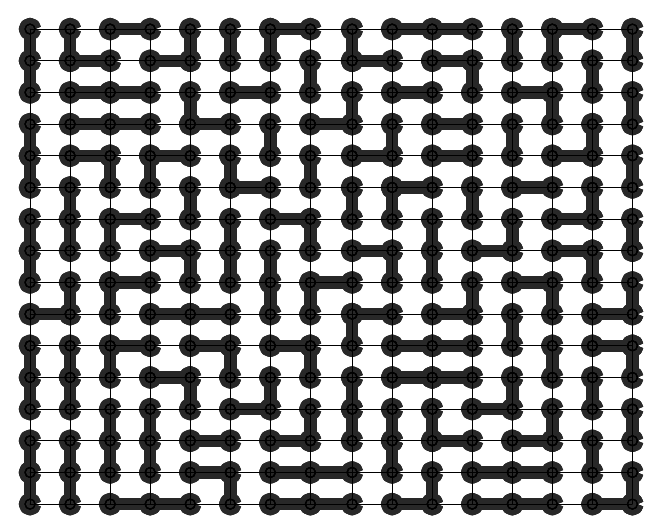}}
\end{center}
\caption{Before (left) and after (right) the Subgraph Reshaping Algorithm 
applied on partitioning given by Parallel Aggregation Algorithm.}
\label{t5-then-rank1}
\end{figure}

Another inspiring observation is that, 
in Table \ref{square-dirichlet}, 
even if we set a limit $t=5$ for the maximal number of vertices in an aggregate, 
the resulting aggregates have an average number of vertices  
ranging from $2.97$ to $2.99$. 
We plot the aggregates of an unweighted graph corresponding to a $16\times 16$ square grid 
on the left of Figure \ref{t5-then-rank1}, 
and observe that, 
some aggregates contain 5 vertices and some contains only 1. 
We then use the rank one optimization discussed in Section \ref{ssec:rank1} 
and apply subgraph reshaping algorithm (Algorithm \ref{alg:sra})  
as a post-process of the GPU parallel aggregation algorithm (Algorithm \ref{alg:paa}),
and plot the resulting aggregates on the right of Figure \ref{t5-then-rank1}. 
Since the subgraph reshaping algorithm does not change the number of aggregates, 
the coarsening ratios on the left and right of Figure \ref{t5-then-rank1} are identical and 
are equal to $2.72$. 
The energy of the $\ell^{2}$ projection is deceased 
from 
$|Q|_{\widehat A}^{2}=2.51$ (left of Figure \ref{t5-then-rank1}) 
to 
$|Q|_{\widehat A}^{2}=2.19$ (right of Figure \ref{t5-then-rank1}). 
However, two level convergence rate increases 
from 
$|E|_{\widehat A}=0.67$ (left of Figure \ref{t5-then-rank1}) 
to 
$|E|_{\widehat A}=0.69$ (right of Figure \ref{t5-then-rank1}). 

  Some more comments on the reshaping algorithm are in order.  For
  isotropic problems, the reshaping does not have significant impact
  of on the convergence rate because aggregation obtained by standard
  approach already results in a good convergence rate.  However, for
  anisotropic problems reshaping improves the convergence rate.  In
  this case, starting with aggregates of arbitrary shape, the
  reshaping procedure results in aggregates aligned with the
  anisotropy and definitely improves the overall convergence rate.

  In addition, even for isotropic case, the numerical results in the
  manuscript indicate that subgraph reshaping can be essential for
  variety of cycling algorithms when aggressive coarsening is applied.
  As shown in Example~\ref{example-1} and Example~\ref{example-2}, the
  aggregation reshaping helps for some isotropic and anisotropic
  problems when coarsening ratio is 8.  In
  Table~\ref{square-dirichlet-allqa} and
  Table~\ref{square-neumann-allqa}, we observe that such coarsening
  ratio can be achieved by skipping every other level in our current
  multilevel hierarchy.

  Clearly, further investigation about the reshaping is needed for more general
  problems that have both anisotropic and isotropic regions.
  Analyzing such cases as well as testing how much improvement in
  the convergence can be achieved by subgraph reshaping for specific
  coarsening and cycling strategies are subject of an ongoing and
  future research.

\subsection{Numerical Tests for GPU Implementation}
In this section, we perform numerical experiments to demonstrate the
efficiency of our proposed AMG method and discuss the specifics
related to the use of GPUs as main platform for computations.  We test
the parallel algorithm on Laplace equation discretized on
quasi-uniform grids in 2D. Our test and comparison platform is the
NVIDIA Tesla C2070 together with a Dell computing workstation. Details
in regard to the machine are given in Table \ref{tab:test_platform}.

\begin{table}[h]
\begin{center}
\begin{tabular}{r|l}
\hline \hline
CPU Type & Intel  \\
CPU Clock & 2.4 GHz \\
Host Memory & 16 GB  \\
\hline 
GPU Type & NVIDIA Tesla C2070 \\
GPU Clock &  575MHz \\
Device Memory & 6 GB \\
CUDA Capability &   2.0\\
\hline
Operating System & RedHat  \\
CUDA Driver &   CUDA 4.1\\
Host Complier &   gcc 4.1\\
Device Complier &  nvcc 4.1\\
CUSP &  v0.3.0\\
\hline \hline
\end{tabular}
\end{center}
\caption{Test Platform}
\label{tab:test_platform}
\end{table}

Because our aim is to demonstrate the improvement of our algorithm
on GPUs, we concentrate on comparing the method we describe here with the
parallel smoothed aggregation AMG method implemented in the CUSP
package \cite{Garland2010}.  

We consider the standard linear finite element method for the Laplace
equation on unstructured meshes.  The results are shown in the
Table~\ref{tab:whole-compare}.  Here, CUSP uses smoothed aggregation
AMG method with V-cycles, and our method is UA-AMG with K-cycles. The
stoping criterion is that 
the $\ell^{2}$ norm of the relative residual is less
than $10^{-6}$.  According to the results, we can see that our
parallel UA-AMG method converges uniformly with respect to the problem
size.  This is due the improved aggregation algorithm constructed
by our parallel aggregation method and the K-cycle used in the solver
phase.  We can see that our method is about 3 to 4 times
faster in setup phase, 
which demonstrate the efficiency of our parallel aggregation
algorithm.  In the solver phase, due to the factor that we use
K-cycle, which does much more work on the coarse grids, our solver
phase is a little bit slower than the solver phase implemented in
CUSP.  However, the use of a K-cycle yields a uniformly convergent
UA-AMG method, which is an essential property for designing scalable
solvers.  When the size of the problem gets larger, we expect the
computational time of our AMG method to scale linearly whereas the AMG
method in CUSP seems to grows faster than linear and will be slower
than our solver phase eventually.  Overall, our new AMG solver is
about $1.5$ times faster than the smoothed aggregation AMG method in
CUSP in terms of total computational time, and the numerical tests
suggests that it converges uniformly for the Poisson problem.

\begin{center}
\begin{table}
\captionsetup{width=28em}
\begin{center}
\begin{tabular}{|c|c|c|c|c|c|c|c|c|}
\hline
& \multicolumn{4}{|c|}{\#DoF = 1 million}&\multicolumn{4}{|c|}{\#DoF = 4 million}\\  \hline
& \# Iter.& Setup & Solve&Total& \# Iter.&Setup&Solve&Total \\ \hline
CUSP & 36& 0.63&0.35&0.98& 41&2.38&1.60&3.98 \\ \hline
New & 19 &0.13&0.47&0.60 & 19 &0.62&2.01&2.63 \\ \hline
\end{tabular}
\caption{Comparison between the parallel AMG method in CUSP package (smoothed aggregation AMG with V-cycles) and our new parallel AMG method (UA-AMG with K-cycles). }\label{tab:whole-compare}
\end{center}
\end{table}
\end{center}

\end{document}